\documentstyle[12pt]{article}

\newcommand{\braid}[3]{{#1}$\lower4pt\hbox{$\oo\atop\raise4pt
           \hbox{$\scriptscriptstyle {#3} $}$}${#2}}
\newcommand{\twist}[3]{{#1}${\,\scriptscriptstyle {#3}}\atop\raise9pt
           \hbox{$\scriptstyle\oo$} ${#2}}

\newcommand{\be}{\begin{eqnarray}}
\newcommand{\ee}{\end{eqnarray}}

\newcommand{\oo}{\otimes  }

\newcommand{\si}{\sigma}

\newcommand{\om}{\omega}

\begin{document}

\begin{titlepage}
\begin{center}
{\Large\bf  Super-jordanian deformation of 
the orthosymplectic Lie superalgebras}
\end{center}
\vskip 1cm
\centerline{\large \bf
P.P. Kulish\footnote{On leave of absence from the 
St.Petersburg Department of the Steklov
Mathematical Institute, Fontanka 27, St.Petersburg, 191011,
Russia.}
\footnote{Partially supported by the RFFI grant N 98-01-00310.}}
\vskip 1mm
\begin{center}
{Dipartimento di Fisica del'Universit\`a and INFN, Firenze}\\
{L.go Fermi 2 - I50125 Firenze, Italy}
\end{center}

\vskip 5mm
\vskip 5mm

\begin{abstract}
The recently proposed jordanian quantization 
of the Lie superalgebra $osp(1|2)$ 
due to the embedding $sl(2) \subset osp(1|2)$,  
is extended including odd generators into the twisting element 
$\cal F$. This deformation is obtained as a contraction 
of the quantum superalgebra ${\cal U}_{q}(osp(1|2))$. 

\end{abstract}

\end{titlepage}


Connection between the Drinfeld-Jimbo deformation \cite{DR,JIM} and 
the jordanian deformation of $sl(2)$ \cite{MAN,ZAK,GER,OG} was pointed 
out already in \cite{GER}. Technically it is given by a similarity 
transformation of the $R$-matrix and $T$-matrix with the element 
$\exp (\theta \, ad\,X_+)$ in the fundamental representation of $sl(2)$  
and contraction $q \to 1$ with fixed $(1-q)\theta = \xi$.  
In particular, one arrives to the jordanian plane 
$x^{\prime }y^{\prime } = y^{\prime }x^{\prime }+\xi {y^{\prime }}^2$ 
of ${\cal U}_{\xi}(sl(2))$ after the transformation: 
$$ 
x^{\prime } = x + \theta y \,, \quad y^{\prime } = y 
$$  
of the $q$-plane $x y = q y x$ \cite{OG}. 

The recently obtained jordanian twist of the universal enveloping 
superalgebra ${\cal U}(osp(1|2))$ \cite{CK} 
was constructed due to the embedding $ sl(2) \subset osp(1|2)$.  
The mentioned above procedure being applied to the $9 \times 9$ 
$R$-matrix of the quantum superalgebra ${\cal U}_q(osp(1|2))$ \cite{KR}
gives extra contribution into 
the twisting element ${\cal F}$ due to the odd generator $v_+$. 
This twisting element is factorised 
$$ 
{\cal F}^{(sj)} = {\cal F}^{(s)}{\cal F}^{(j)} \,. 
$$  
The jordanian twisting element \cite{GER,OG} 
\be
\label{jt} 
{\cal F}^{(j)} = \exp \left(\frac 12 h \oo \ln(1 + 2 \xi X_+)\right) = 
e^{\si}  
\ee
depends on the Borel subalgebra generators 
$h,\, X_+ \in B_+ \subset sl(2) \subset osp(1|2)$ 
while the super-part ${\cal F}^{(s)}$ depends on $v_+$ only 
(with $X_+ = 4 v_+^2$ and $[h,\,X_+] = 2 X_+$).

Recall the twisting of Hopf algebras. 
A Hopf algebra ${\cal A}(m, \Delta, \epsilon, S)$ with
multiplication $m\colon {\cal A}\otimes {\cal A}\to {\cal A}$,
coproduct 
$\Delta \colon {\cal A}\to {\cal A} \otimes {\cal A}$, 
counit $\epsilon \colon {\cal A}\to C$, 
and antipode $S\colon {\cal A}\to {\cal A}$ 
(see definitions in Refs.\cite{DR,FAD})
can be transformed with an invertible element ${\cal F} \in 
{\cal A} \otimes {\cal A}$, ${\cal F}=\sum f_i^{(1)}\otimes f_i^{(2)}$ 
into a twisted one ${\cal A}_t(m,\Delta _t,\epsilon ,S_t)$ \cite{DR2}.
This Hopf algebra ${\cal A}_t$ has the same multiplication 
and counit maps but the twisted coproduct and antipode
$$
\Delta _t(a)={\cal F}\Delta (a){\cal F}^{-1},\quad S_t(a)=vS(a)v^{-1},\quad
v=\sum f_i^{(1)}S(f_i^{(2)}),\quad a\in {\cal A}.
$$
The twisting element has to satisfy the equations 
\begin{eqnarray}(\epsilon \otimes  id)({\cal F}) = 
(id \otimes  \epsilon)({\cal F})=1, 
\end{eqnarray} 
\begin{eqnarray}
{\cal F}_{12}(\Delta \otimes  id)({\cal F}) =
{\cal F}_{23}(id \otimes  \Delta)({\cal F})\,. 
\label{TE}
\end{eqnarray}

A quasitriangular Hopf algebra 
${\cal A}(m,\Delta , \epsilon , S,{\cal R})$ has
additionally an element ${\cal R}\in {\cal A}\otimes {\cal A}$ 
(a universal $R$-matrix) \cite{DR}, which relates the coproduct $\Delta $ 
and its opposite coproduct $\Delta^{op} $ by the similarity transformation 
$$ 
\Delta^{op}(a) = {\cal R} \Delta(a){\cal R}^{-1}\,, \quad a \in {\cal A}\,.  
$$ 
A twisted quasitriangular quantum
algebra ${\cal A}_t(m, \Delta _t, \epsilon, S_t, {\cal R}_t)$ 
has the twisted universal $R$-matrix
\begin{eqnarray}{\cal R}_t=\tau  ({\cal F})\,{\cal R}\,{\cal F}^{-1},
\label{Rt}
\end{eqnarray}where $\tau $ means the permutation of the tensor factors: 
$\tau (f\otimes g)=(g\otimes f)\,, \tau ({\cal F})={\cal F}_{21}$.

Note that the composition of appropriate twists can be defined
${\cal F} = {\cal F}_2 {\cal F}_1$. The element ${\cal F}_1$ has to
satisfy the twist equation with the coproduct of the original Hopf algebra,
while ${\cal F}_2$ must be its solution for $\Delta_{t_1}$ of the
intermediate Hopf algebra twisted by ${\cal F}_1$.
In particular, if ${\cal F}$ is a solution to the twist equation (\ref{TE})
then ${\cal F}^{-1}$ satisfies this equation with $\Delta \rightarrow \Delta
_t$.

The quasitriangular Hopf superalgebra ${\cal U}_q(osp(1|2))$ \cite{KR} 
is generated by three elements $\{h, v_-, v_+\}$,  
analogously to the universal enveloping of $osp(1|2)$ or $sl(2)$, 
subject to the relations 
\begin{eqnarray}
\label{q-osp}
[ h, v_{\pm} ] = \, \pm \; v_{\pm} ,\hskip 1.5cm [ v_+, v_-] = \, - \;
\frac{1}{4} \,(q^h - q^{-h})/ (q - q^{-1})  \,,  
\end{eqnarray}
\noindent 
where the commutator $[\,,\,]$ is understood as the $Z_2$-graded one: 
$$ 
[a, b] = ab - (-1)^{p(a)p(b)}ba \,, 
$$ 
with $p(a) = 0, 1$ being the parity (even or odd) of the element. The 
coproduct is 
\begin{eqnarray}
\Delta_{q} (h) &=& h \otimes 1 + 1 \otimes h \;, \nonumber \\
\Delta_{q} (v_{\pm}) &=& v_{\pm} \otimes q^{h/2} + 
q^{-h/2} \otimes v_{\pm} \;.  
\label{qco} 
\end{eqnarray} 
The $osp(1|2)$ commutation relations follow from (\ref{q-osp}) in the 
limit $q \to 1$ and adding $X_{\pm} = \pm 4 (v_{\pm})^2$ as 
the Lie superalgebra generators. 
It is worthy to note that, while $sl(2)$ is embedded into $osp(1|2)$,
such embedding does not exist for $sl_q(2)$ into $osp_q(1|2)$
because the coproduct of even elements $X_{\pm} \sim v_{\pm}^2$ 
$$ 
\Delta_q (X_{\pm}) \sim X_{\pm} \oo q^h + q^{-h} \oo X_{\pm} + 
(q^{\pm 2} - 1)q^{-h}v_{\pm} \oo v_{\pm}q^h 
$$ 
includes also odd generators.

The jordanian twist ${\cal F}^{(j)}$ preserving the algebraic relations 
among the generators of  ${\cal U}(osp(1|2))$, results in the twisted 
coproduct $\Delta_{j}$ \cite{CK}  
\begin{eqnarray}
\Delta_{j} (h) &=& h \otimes e^{-2\sigma} + 1 \otimes h \;, \nonumber \\
\Delta_{j} (v_+) &=& v_+ \otimes e^{\sigma} + 1 \otimes v_+ \;, \label{def} \\
\Delta_j (v_-) &=& v_- \otimes e^{-\sigma} 
     + 1 \otimes v_- +  \xi\, h \otimes v_+ e^{-2\sigma} .  \nonumber
\end{eqnarray}

The $R$-matrix formulation of the quantum algebra ${\cal U}_q(osp(1|2))$ 
is related with the $9 \times 9$ matrix \cite{KR}:  
\begin{equation}
\label{R-mat}R=\left(
\begin{array}{ccccccccc}
q & 0 & 0 & 0 & 0 & 0 & 0 & 0 & 0 \\
0 & 1 & 0 & a & 0 & 0 & 0 & 0 & 0\\
0 & 0 & q^{-1} & 0 & b & 0 & e & 0 & 0\\
0 & 0 & 0 & 1 & 0 & 0 & 0 & 0 & 0\\
0 & 0 & 0 & 0 & 1 & 0 & c & 0 & 0\\
0 & 0 & 0 & 0 & 0 & 1 & 0 & d & 0\\
0 & 0 & 0 & 0 & 0 & 0 & q^{-1} & 0 & 0\\
0 & 0 & 0 & 0 & 0 & 0 & 0 & 1 & 0\\
0 & 0 & 0 & 0 & 0 & 0 & 0 & 0 & q
\end{array}
\right) ,
\end{equation}
where non zero entries are 
$$ 
a = d = \om = q - q^{-1} \,, \quad b = c = - \om / \sqrt{q} \,, \quad 
e = \om (1+1/q)\,.  
$$ 
Following the $sl(N)$ case \cite{GER,KLM}, let us do 
the similarity transformation of this $R$-matrix with 
the tensor square $M \oo M$ 
of the matrix $M = I + \theta \rho (X_+)$ 
\begin{equation}
\label{M-matr}M=\left(
\begin{array}{ccc}
1 & 0 & \theta \\
0 & 1 & 0 \\
0 & 0 & 1
\end{array}
\right) \,. 
\end{equation}
The transformed $R$-matrix satisfying the graded 
Yang-Baxter equation (YBE) \cite{KSK}, 
gets nine more non zero entries: 
\begin{equation}
\label{tR-mat}R=\left(
\begin{array}{ccccccccc}
q & 0 & x_1 & 0 & x_2 & 0 & x_3 & 0 & x_4 \\
0 & 1 & 0 & a & 0 & x_5 & 0 & 0 & 0\\
0 & 0 & q^{-1} & 0 & b & 0 & e & 0 & x_6\\
0 & 0 & 0 & 1 & 0 & 0 & 0 & x_7 & 0\\
0 & 0 & 0 & 0 & 1 & 0 & c & 0 & x_8\\
0 & 0 & 0 & 0 & 0 & 1 & 0 & d & 0\\
0 & 0 & 0 & 0 & 0 & 0 & q^{-1} & 0 & x_9\\
0 & 0 & 0 & 0 & 0 & 0 & 0 & 1 & 0\\
0 & 0 & 0 & 0 & 0 & 0 & 0 & 0 & q
\end{array}
\right) ,
\end{equation} 
which are the following: 
$$ 
\begin{array}{lllllllll}
x_1 &=& -\om\, \theta \,, & x_2 &=& b\, \theta \,, & x_3 &=& 
\om\, \theta/q \,, \\  
x_4 &=& (\om\, \theta)^2/(1+q) \,, & 
x_5 &=& -\om\, \theta \,, & x_6 &=& -\om\, \theta/q \,, \\
x_7 &=& \om\, \theta \,, &  x_8 &=& -c\, \theta \,, & x_9 &=& 
\om\, \theta \,.  
\end{array}
$$ 
All the new entries are proportional to $ \om \, \theta $. 
Hence one can consider a limit: $\theta = \xi / \om \,,\, q \to 1$. 
The limiting (super-jordanian) $R$-matrix $R^{(sj)}(\xi)$ is also 
a solution to the graded YBE. The form of $R^{(sj)}(\xi)$ 
is different from the jordanian $R$-matrix 
which was used for a triangular twist of $osp(1|2)$   
\begin{equation}
\label{sjR-mat}R^{(sj)}(\xi) \, =\left(
\begin{array}{ccccccccc}
1 & 0 & -\xi & 0 & -\xi & 0 & \xi & 0 & \xi^2/2 \\
0 & 1 & 0 & 0 & 0 & -\xi & 0 & 0 & 0\\
0 & 0 & 1 & 0 & 0 & 0 & 0 & 0 & -\xi\\
0 & 0 & 0 & 1 & 0 & 0 & 0 & \xi & 0\\
0 & 0 & 0 & 0 & 1 & 0 & 0 & 0 & \xi\\
0 & 0 & 0 & 0 & 0 & 1 & 0 & 0 & 0\\
0 & 0 & 0 & 0 & 0 & 0 & 1 & 0 & \xi\\
0 & 0 & 0 & 0 & 0 & 0 & 0 & 1 & 0\\
0 & 0 & 0 & 0 & 0 & 0 & 0 & 0 & 1
\end{array}
\right) \,. 
\end{equation}

One can now take the upper triangular $L^{(+)}$ 
matrix of the FRT-formalism \cite{FAD}  
\begin{equation}
\label{L-mat}L^{(+)}=\left(
\begin{array}{ccc}
E^{-1} & V & H \\
0 & 1 & W \\
0 & 0 & E
\end{array}
\right) = (\rho \oo id) {\cal R}^{(sj)} \,,  
\end{equation}
and define the commutation relations of the generators $H, E, V, W$ 
from the FRT-relation 
$$ 
R^{(sj)}(\xi) L^{(+)}_1 L^{(+)}_2 = L^{(+)}_2 L^{(+)}_1 R^{(sj)}(\xi) \,,  
$$ 
where the $Z_2$-graded tensor product is used defining 
$L^{(+)}_1 = L^{(+)} \oo I$ and $ L^{(+)}_2 = I \oo L^{(-)} $ \cite{KSK}. 
These commutation relations are ($V$ and $W$ are odd) 
$$ 
[E, V] = 0 \,, \quad [E, W] = 0 \,, 
$$
$$ 
[H, E] = \xi (E^2 - 1)\,,\quad [H, V] = \xi (V (E^{-1}- E) - W)\,, 
$$
$$ 
[V, V] = \xi (1 - E^{-2})\,,\quad [W, W] = \xi (E^2 - 1)\,, 
$$
$$ 
VW + WV = - \xi (E - E^{-1})\,,\quad 
V^2 + W^2 = \frac 12 \xi (E^2 - E^{-2})\,. 
$$ 
It is easy to see that one can take $H$ and $W$ as two independent 
generators, while $V = - W E^{-1}$ and $E^2 = 1 + 2W^2/\xi$. 

The super-jordanian $R$-matrix is triangular 
$R_{21}^{(sj)}(\xi)\,R^{(sj)}(\xi) = 1$ and it would be interesting 
to find the corresponding twisting element. Let us start from 
the twist matrix $F = (\rho \oo \rho) {\cal F}$ in the fundamental 
representation $\rho$. Having in mind embedding of the jordanian 
twist of $sl(2)$ into $osp(1|2)$ \cite{CK} and 
a possibility of twist composition, 
we are looking for the super-jordanian 
twist as the product ${\cal F}^{(sj)} = {\cal F}^{(s)}{\cal F}^{(j)}$ 
with corresponding matrices: block diagonal 
\begin{equation}
\label{jF-mat}F^{(j)}= diag \left( \left( 
\begin{array}{ccc} 
1 & 0 & \xi  \\
0 & 1 & 0 \\
0 & 0 & 1 
\end{array}
\right),\, \left(
\begin{array}{ccc} 
1 & 0 & 0  \\
0 & 1 & 0 \\
0 & 0 & 1 
\end{array}
\right),\, \left( 
\begin{array}{ccc} 
1 & 0 & -\xi  \\
0 & 1 & 0 \\
0 & 0 & 1 
\end{array}
\right) 
\right) = (\rho \oo \rho) e^{h \oo \sigma} ,  
\end{equation} 
for the jordanian twist, and 
\begin{equation}
\label{sF-mat}F^{(s)}=\left(
\begin{array}{ccccccccc}
1 & 0 & 0 & 0 & \xi/2 & 0 & 0 & 0 & -\xi^2/8 \\
0 & 1 & 0 & 0 & 0 & \xi/2 & 0 & 0 & 0\\
0 & 0 & 1 & 0 & 0 & 0 & 0 & 0 & 0\\
0 & 0 & 0 & 1 & 0 & 0 & 0 & -\xi/2 & 0\\
0 & 0 & 0 & 0 & 1 & 0 & 0 & 0 & -\xi/2\\
0 & 0 & 0 & 0 & 0 & 1 & 0 & 0 & 0\\
0 & 0 & 0 & 0 & 0 & 0 & 1 & 0 & 0\\
0 & 0 & 0 & 0 & 0 & 0 & 0 & 1 & 0\\
0 & 0 & 0 & 0 & 0 & 0 & 0 & 0 & 1
\end{array}
\right) = (\rho \oo \rho) {\cal F}^{(s)} \,, 
\end{equation} 
for the super-part which depends on $v$ only and has 
a property: $F_{21}^{(s)} = (F^{(s)})^{-1} $. 
It is easy to check that 
$$
R^{(sj)}(\xi) = 
F_{21}^{(s)}\, F_{21}^{(j)}\, (F^{(j)})^{-1}\, 
(F^{(s)})^{-1}\,. 
$$

Taking into account the property of the universal $R$-matrix $\cal R$ 
\cite{DR}, and the definition 
$L^{(+)} = (\rho \oo id) {\cal R}^{(sj)}$, one gets the coproducts of the 
$L^{(+)}$ entries 
$$ 
\begin{array}{lll}
\Delta(E) &=& E \oo E\,, \\ 
\Delta(V) &=& V \oo E^{-1} + 1 \oo V \,, \\ 
\Delta(W) &=& W \oo 1 + E \oo W \,, \\ 
\Delta(H) &=& H \oo E^{-1} + E \oo H - W \oo V \,.  
\end{array}
$$  

We have to find the expressions of the generators $H, E, V, W$ in terms 
of $h, v$, and define the super-twist ${\cal F}^{(s)}$ from the 
intertwining relation. The form of this element ${\cal F}^{(s)}$ in 
question is 
$$ 
{\cal F}^{(s)} = 
\exp \left(-2 \xi\, (v \oo v)\, 
\varphi (\si \oo 1\,, 1 \oo \si)\right) \,, 
$$ 
where $\varphi (\si_1,\, \si_2)$ is a symmetric function of its 
arguments. Fixing this structure of ${\cal F}^{(s)}$ and considering 
the fundamental representation $\rho$ only for the first factor 
in ${\cal A} \oo {\cal A}$, 
we get the following expressions of  the generators $H, E, V, W$ 
\begin{equation}
\label{t-gen} 
H = \xi h e^{\si} - 2(\xi v)^2 e^{-\si}, \; 
E = \exp ( \frac 12 \ln (1 + 2\xi X) ) = e^{\sigma}, \; 
V = -2\xi v e^{-\si} , \; W = 2\xi v \,.  
\end{equation}
Hence the coproducts of the generators $h, v$ are the following 
\begin{eqnarray}
\Delta_{sj} (h) &=& h \oo E^{-2} + 1 \oo h + 
4 \xi\, v E^{-1} \oo v E^{-2}\,, \nonumber \\   
\Delta_{sj} (v) &=& v \oo 1 + e^{\si} \oo v \,. \label{Dsj} 
\end{eqnarray} 
Thus the super-jordanian twist of the Borel subalgebra $sB_+$ 
of $osp(1|2)$ is defined. 

To define the corresponding deformation of the ${\cal U}(osp(1|2))$  
we have to find ${\cal F}^{(sj)}$  and the coproduct 
of the generator $v_-$ 
\be 
{\cal F}^{(sj)} \Delta (v_-) ({\cal F}^{(sj)})^{-1} = 
{\cal F}^{(s)} (v_- \oo E^{-1} + 1 \oo v_- + \xi\, h \oo v_+ E^{-2}) 
({\cal F}^{(s)})^{-1}\,. 
\ee 
The knowledge of the exact form of the super-part of the twisting 
element seems necessary to get a final expression for 
$\Delta_{sj}(v_-)$. Although one can prove according to \cite{FR} 
(with appropriate modifications due to the $Z_2$-grading) 
the existence of ${\cal F}^{(s)}$ order by order in $\xi$ 
we do not have a closed form. The intertwining relation 
$$
{\cal F}^{(s)} \Delta_{j} (v_+) ({\cal F}^{(s)})^{-1} = 
(v_+ \oo E + 1 \oo v_+) 
({\cal F}^{(s)})^{-2} = v_+ \oo 1 + E \oo v_+
$$ 
also can be used to find ${\cal F}^{(s)}$. The conjectured form of 
$\varphi (\si_1,\, \si_2)$ is the following 
$$ 
\varphi (\si \oo 1\,, 1 \oo \si) =  \sum_{k=1}^{\infty} f_k(\si) 
\oo f_k(\si) \,. 
$$ 
Each $f_k(\si)$  of this expression is characterised by its non zero 
contribution starting from the irreducible representation of 
spin $s = k/2,\, dim V_s = 4s + 1$ and by the first term $(\xi X_+)^{k-1}$ 
with a coefficient to be defined. In particular, one gets 
$f_1(\si) = 2/(e^{\si} + 1)$. 

Taking the generators of the Lie superalgebra as invariant vector 
fields on the group supermanifold, the twist ${\cal F}$ 
as a bidifferential operator permits to put forward 
the deformation quantization approach of \cite{FLA}. 
Using the explicit form of the $9 \times 9$ $R$-matrix (\ref{sjR-mat}), 
one can also study the dual Hopf superalgebra (a quantum supergroup) 
according to the relation 
$$ 
R^{(sj)}(\xi) T_1 T_2 = T_2 T_1 R^{(sj)}(\xi). 
$$ 
Corresponding formulas, 
together with the classical $r$-matrix, are given in \cite{JS}. 
One can apply the deformed $osp(1|2)$ to the oscillator like realizations 
of quantum superalgebras \cite{CPT}, and to integrable models related to 
the orthosymplectic superalgebras \cite{K}. 

${\bf Acknowledgements}$. 
The author is grateful to E. Celeghini for useful 
discussions on superalgebras. 
The hospitality and support of the Institute for 
Theoretical Physics of the Wroclaw University, 
the Physics Department of the Florence University 
and the Florence section of INFN are acknowledged.

\end{document}